\pdfoutput=1
\RequirePackage{ifpdf}
\ifpdf 
\documentclass[pdftex]{sigma}
\else
\documentclass{sigma}
\fi

\usepackage{bm}

\numberwithin{equation}{section}
\newtheorem{Theorem}{Theorem}[section]
\newtheorem{Corollary}[Theorem]{Corollary}
\newtheorem{Lemma}[Theorem]{Lemma}
\newtheorem{Proposition}[Theorem]{Proposition}
{\theoremstyle{definition}
\newtheorem{Definition}[Theorem]{Definition}
\newtheorem{Remark}[Theorem]{Remark}
}

\begin{document}

\allowdisplaybreaks

\renewcommand{\thefootnote}{$\star$}

\renewcommand{\PaperNumber}{021}

\FirstPageHeading

\ShortArticleName{Commutative Families of the Elliptic Macdonald Operator}

\ArticleName{Commutative Families\\ of the Elliptic Macdonald Operator\footnote{This paper is a~contribution
to the Special Issue in honor of Anatol Kirillov and Tetsuji Miwa.
The full collection is available at \href{http://www.emis.de/journals/SIGMA/InfiniteAnalysis2013.html}
{http://www.emis.de/journals/SIGMA/InfiniteAnalysis2013.html}}}

\Author{Yosuke SAITO}

\AuthorNameForHeading{Y.~Saito}

\Address{Mathematical Institute of Tohoku University, Sendai, Japan}
\Email{\href{yousukesaitou7@gmail.com}{yousukesaitou7@gmail.com}}

\ArticleDates{Received October 01, 2013, in f\/inal form February 25, 2014; Published online March 11, 2014}

\Abstract{In the paper~[\textit{J.~Math. Phys.} \textbf{50} (2009), 095215, 42~pages],
Feigin, Hashizume, Hoshino, Shiraishi, and Yanagida constructed two families of commuting operators which contain the Macdonald operator (commutative families of the Macdonald operator). They used the Ding--Iohara--Miki algebra and the trigonometric Feigin--Odesskii algebra.
In the previous paper~[arXiv:1301.4912], the present author constructed the elliptic Ding--Iohara--Miki al\-gebra and the free f\/ield
realization of the elliptic Macdonald operator.
In this paper, we show that by using the elliptic Ding--Iohara--Miki algebra and the elliptic Feigin--Odesskii algebra,
we can construct commutative families of the elliptic Macdonald operator.
In Appendix, we will show a~relation between the elliptic Macdonald operator and its kernel function by the free f\/ield
realization.}


\Classification{17B37; 33D52}

\renewcommand{\thefootnote}{\arabic{footnote}}
\setcounter{footnote}{0}

\textbf{Notations.} In this paper, we use the following symbols.
\begin{gather*}
\mathbf{Z}: \text{the set of integers},
\quad
\mathbf{Z}_{\geq{0}}:=\{0, 1, 2,\dots\},
\quad
\mathbf{Z}_{>0}:=\{1, 2,\dots\},
\\
\mathbf{Q}: \text{the set of rational numbers},
\quad
\mathbf{Q}(q,t): \text{the f\/ield of rational functions of $q$, $t$ over $\mathbf{Q}$},
\\
\mathbf{C}: \text{the set of complex numbers},
\quad
\mathbf{C}^{\times}:=\mathbf{C}\setminus\{0\},
\\
\mathbf{C}[[z,z^{-1}]]: \text{the set of formal power series of $z$, $z^{-1}$ over $\mathbf{C}$}.
\end{gather*}

A~partition is a sequence $\lambda=(\lambda_{1},\dots,\lambda_{N})\in(\mathbf{Z}_{\geq{0}})^{N}$, $N\in\mathbf{Z}_{>0}$ satisfying the condition $\lambda_{i}\geq\lambda_{i+1}$ for each~$i$, $1\leq{i}\leq{N-1}$.
We denote the set of partitions by~$\mathcal{P}$.
For a~partition~$\lambda$,
$\ell(\lambda):=\sharp\{i: \lambda_{i}\neq 0\}$ denotes the length of $\lambda$ and
$|\lambda|:=\sum\limits_{i=1}^{\ell(\lambda)}\lambda_{i}$ denotes the size of~$\lambda$.

Let $q,p\in\mathbf{C}$ be complex parameters satisfying $|q|<1$, $|p|<1$.
We def\/ine the $q$-inf\/inite product as $(x;q)_{\infty}:=\prod\limits_{n\geq{0}}(1-xq^{n})$ and the theta function as
\begin{gather*}
\Theta_{p}(x):=(p;p)_{\infty}(x;p)_{\infty}\big(px^{-1};p\big)_{\infty}.
\end{gather*}
We set the double inf\/inite product as $(x;q,p)_{\infty}:=\prod\limits_{m,n\geq{0}}(1-xq^{m}p^{n})$ and the elliptic
gamma function as
\begin{gather*}
\Gamma_{q,p}(x):=\frac{(qpx^{-1};q,p)_{\infty}}{(x;q,p)_{\infty}}.
\end{gather*}

\section{Introduction}

The elliptic Feigin--Odesskii algebra $\mathcal{A}(p)$ was originally
introduced by Feigin and Odesskii~\cite{7} as an algebra generated by a~certain family of multivariate elliptic
functions.
The algebra $\mathcal{A}(p)$ has a~product called the star product~$\ast$ (Def\/inition~\ref{def2.7}), and in fact
$\mathcal{A}(p)$ is unital, associative, and commutative in terms of the star product~$\ast$.

On the other hand, Miki~\cite{6} constructed a~$q$-deformation of the $W_{1+\infty}$ algebra which has a~structure of
Ding and Iohara's quantum group~\cite{5}.
Miki's quantum group is also obtained from the free f\/ield realization of the Macdonald operator~\cite{1}.
We call the quantum group the Ding--Iohara--Miki algebra $\mathcal{U}(q,t)$ (Def\/inition~\ref{def1.1}).
Then Feigin, Hashizume, Hoshino, Shiraishi, and Yanagida showed the following~\cite{1}:
\begin{itemize}\itemsep=-2pt
\item
There exists a~trigonometric degeneration of the elliptic Feigin--Odesskii algebra $\mathcal{A}(p)$.
We call the algebra the trigonometric Feigin--Odesskii algebra $\mathcal{A}$ (Def\/inition~\ref{def1.6}).
\item
By using the Ding--Iohara--Miki algebra $\mathcal{U}(q,t)$ and the trigonometric Feigin--Odesskii algebra
$\mathcal{A}$, we can obtain commutative families of the Macdonald operator.
\end{itemize}
The trigonometric Feigin--Odesskii algebra $\mathcal{A}$ also has the star product~$\ast$.
In~\cite{1}, the commutativity of~$\mathcal{A}$ with respect to the star product $\ast$ is shown by using some combinatorial tools such as the Gordon f\/iltrations.
These tools are also available in the elliptic case.
Therefore we can say a combinatorial way to prove the commutativities of both trigonometric and elliptic Feigin--Odesskii algebras has been found in~\cite{1}.

The aim of this paper is to extend some results of Feigin, Hashizume, Hoshino, Shiarishi, and Yanagida~\cite{1} to the elliptic case. That is, we construct commutative families of the elliptic Macdonald operator by using the elliptic Ding--Iohara--Miki algebra $\mathcal{U}(q,t,p)$ and the elliptic Feigin--Odesskii algebra~$\mathcal{A}(p)$. We utilize our previous result~\cite{3}, the free f\/ield realization of the elliptic Macdonald operator and the elliptic Ding--Iohara--Miki algebra $\mathcal{U}(q,t,p)$, for this purpose.

In the trigonometric case, due to the theory of the Macdonald symmetric functions, more properties about the
trigonometric Feigin--Odesskii algebra $\mathcal{A}$ and commutative families of the Macdonald operator have been known.
Details can be found in~\cite{1,8}.

\textbf{Organization of this paper.} In Section~\ref{Sec1}, we review the trigonometric case treated in~\cite{1}.
In Section~\ref{Sec2}, f\/irst we recall related materials of the elliptic Ding--Iohara--Miki algebra and the free f\/ield
realization of the elliptic Macdonald operators.
Then using the elliptic Ding--Iohara--Miki algebra and the elliptic Feigin--Odesskii algebra, we derive commutative families of the elliptic Macdonald operators.

In Appendix, we show a~functional equation of the elliptic kernel function with the aid of the free f\/ield realization of the elliptic Macdonald operators.

\section{Review of the trigonometric case}\label{Sec1}

In this section, we review the construction of the commutative families of the Macdonald
operator by Feigin, Hashizume, Hoshino, Shiraishi, and Yanagida~\cite{1}.

\subsection[Ding--Iohara--Miki algebra $\mathcal{U}(q,t)$]{Ding--Iohara--Miki algebra $\boldsymbol{\mathcal{U}(q,t)}$}

The Ding--Iohara--Miki algebra is a~quantum group obtained
from the free f\/ield realization of the Macdonald operator~\cite{1}.
Let $q,t\in\mathbf{C}$ be parameters satisfying $|q|<1$.

\begin{Definition}
[Ding--Iohara--Miki algebra $\mathcal{U}(q,t)$]
\label{def1.1}
Let us def\/ine the structure function $g(x)\in\mathbf{C}[[x]]$ as
\begin{gather*}
g(x):=\frac{(1-qx)(1-t^{-1}x)(1-q^{-1}tx)}{(1-q^{-1}x)(1-tx)(1-qt^{-1}x)}.
\end{gather*}
Let $C$ be a~central, invertible element and $x^{\pm}(z):=\sum\limits_{n\in{\mathbf{Z}}}x^{\pm}_{n}z^{-n}$,
$\psi^{\pm}(z):=\sum\limits_{\pm n\geq{0}}\psi^{\pm}_{n}z^{-n}$ be currents satisfying the relations
\begin{gather*}
[\psi^{\pm}(z), \psi^{\pm}(w)]=0,
\qquad
\psi^{+}(z)\psi^{-}(w)=\frac{g(Cz/w)}{g(C^{-1}z/w)}\psi^{-}(w)\psi^{+}(z),
\\
\psi^{\pm}(z)x^{+}(w)=g\left(C^{\pm\frac{1}{2}}\frac{z}{w}\right)x^{+}(w)\psi^{\pm}(z),
\qquad
\psi^{\pm}(z)x^{-}(w)=g\left(C^{\mp\frac{1}{2}}\frac{z}{w}\right)^{-1}x^{-}(w)\psi^{\pm}(z),
\\
x^{\pm}(z)x^{\pm}(w)=g\left(\frac{z}{w}\right)^{\pm 1}x^{\pm}(w)x^{\pm}(z),
\\
[x^{+}(z),x^{-}(w)]=\frac{(1{-}q)(1{-}t^{-1})}{1{-}qt^{-1}}
\bigg\{\delta\left(C\frac{w}{z}\right)\psi^{+}\big(C^{1/2}w\big)-\delta\left(C^{-1}\frac{w}{z}\right)\psi^{-}\big(C^{-1/2}w\big)\bigg\}.
\end{gather*}
Here we set the delta function by $\delta(x):=\sum\limits_{n\in\mathbf{Z}}x^{n}$.
We def\/ine the Ding--Iohara--Miki algeb\-ra~$\mathcal{U}(q,t)$ to be an associative $\mathbf{C}$-algebra generated by
$\{x^{\pm}_{n}\}_{n\in{\mathbf{Z}}}$, $\{\psi^{\pm}_{n}\}_{n\in{\mathbf{Z}}}$, and $C$.
\end{Definition}

The free f\/ield realization of the Ding--Iohara--Miki algebra is known as follows
In the follo\-wing, let $q,t\in\mathbf{C}$ and we assume $|q|<1$.
First we def\/ine the algebra $\mathcal{B}$ of bosons to be an associative $\mathbf{C}$-algebra generated by
$\{a_{n}\}_{n\in{\mathbf{Z}\setminus\{0\}}}$ and the relation
\begin{gather*}
[a_{m},a_{n}]=m\frac{1-q^{|m|}}{1-t^{|m|}}\delta_{m+n,0}, \qquad m,n\in\mathbf{Z}\setminus\{0\}.
\end{gather*}
We set the normal ordering $\bm{:} \bullet \bm{:}$ as
\begin{gather*}
\bm{:}a_{m}a_{n}\bm{:}=
\begin{cases}
a_{m}a_{n}, & m<n,
\\
a_{n}a_{m}, & m\geq n .
\end{cases}
\end{gather*}
Let $|0 \rangle$ be the vacuum vector which satisf\/ies $a_{n}|0 \rangle=0$, $n>0$.
For a~partition $\lambda$, we set $a_{-\lambda}:=a_{-\lambda_{1}}\cdots a_{-\lambda_{\ell(\lambda)}}$ and def\/ine the
boson Fock space $\mathcal{F}$ as the left $\mathcal{B}$ module
\begin{gather*}
\mathcal{F}:=\operatorname{span}\{a_{-\lambda}|0 \rangle: \lambda\in{\mathcal{P}}\}.
\end{gather*}
Let $\langle 0|$ be the dual vacuum vector which satisf\/ies the condition $\langle 0|a_{n}=0$, $n<0$, and def\/ine the dual
boson Fock space $\mathcal{F}^{\ast}$ as the right $\mathcal{B}$ module
\begin{gather*}
\mathcal{F}^{\ast}:=\operatorname{span}\{\langle 0|a_{\lambda}: \lambda\in\mathcal{P}\},
\qquad
a_{\lambda}:=a_{\lambda_{1}}\cdots a_{\ell(\lambda)}.
\end{gather*}
We def\/ine $n_{\lambda}(a):=\sharp\{i: \lambda_{i}=a\}$ and $z_{\lambda}$, $z_{\lambda}(q,t)$ as
\begin{gather*}
z_{\lambda}:=\prod\limits_{a\geq{1}}a^{n_{\lambda}(a)}n_{\lambda}(a)!,
\qquad
z_{\lambda}(q,t):=z_{\lambda}\prod\limits_{i=1}^{\ell(\lambda)}\frac{1-q^{\lambda_{i}}}{1-t^{\lambda_{i}}}.
\end{gather*}
We def\/ine a~bilinear form $\langle \bullet|\bullet \rangle: \mathcal{F}^{\ast} \times \mathcal{F} \to \mathbf{C}$ by the
following conditions:
\begin{gather*}
(1)
\quad
\langle 0|0 \rangle=1,
\\
(2)
\quad
\langle 0|a_{\lambda}a_{-\mu}|0 \rangle=\delta_{\lambda\mu}z_{\lambda}(q,t).
\end{gather*}

\begin{Proposition}[free f\/ield realization of the Ding--Iohara--Miki algebra $\mathcal{U}(q,t)$]
Set $\gamma:=(qt^{-1})^{-1/2}$ and define
operators $\eta(z)$, $\xi(z)$, $\varphi^{\pm}(z): \mathcal{F} \to \mathcal{F}\otimes\mathbf{C}[[z,z^{-1}]]$ as
\begin{gather*}
\eta(z):=\bm{:}\exp\left(-\sum\limits_{n\neq{0}}(1-t^{n})a_{n}\frac{z^{-n}}{n}\right)\bm{:},
\qquad
\xi(z):=\bm{:}\exp\left(\sum\limits_{n\neq{0}}(1-t^{n})\gamma^{|n|}a_{n}\frac{z^{-n}}{n}\right)\bm{:},
\\
\varphi^{+}(z):=\bm{:}\eta(\gamma^{1/2}z)\xi(\gamma^{-1/2}z)\bm{:},
\qquad
\varphi^{-}(z):=\bm{:}\eta(\gamma^{-1/2}z)\xi(\gamma^{1/2}z)\bm{:}.
\end{gather*}
Then the map
\begin{gather*}
C \mapsto \gamma,
\qquad
x^{+}(z) \mapsto \eta(z),
\qquad
x^{-}(z) \mapsto \xi(z),
\qquad
\psi^{\pm}(z) \mapsto \varphi^{\pm}(z)
\end{gather*}
gives a~representation of the Ding--Iohara--Miki algebra $\mathcal{U}(q,t)$.
\end{Proposition}

Here we collect some notations of symmetric polynomials and symmetric functions~\cite{4}.
Let $q,t\in{\mathbf{C}}$ be parameters and assume $|q|<1$.
We denote the $N$-th symmetric group by $\mathfrak{S}_{N}$ and def\/ine
$\Lambda_{N}(q,t):=\mathbf{Q}(q,t)[x_{1},\dots,x_{N}]^{\mathfrak{S}_{N}}$ as the space of $N$-variables symmetric
polynomials over $\mathbf{Q}(q,t)$.
For $\alpha=(\alpha_{1}, \dots, \alpha_{N})\in{(\mathbf{Z}_{\geq{0}})^{N}}$, we set $x^{\alpha}:=x_{1}^{\alpha_{1}}\cdots
x_{N}^{\alpha_{N}}$.
For a~partition $\lambda$, we def\/ine the monomial symmetric polynomial $m_{\lambda}(x)$ as follows
\begin{gather*}
m_{\lambda}(x):=\sum\limits_{\text{$\alpha: \alpha$ is a~permutation of $\lambda$}}x^{\alpha}.
\end{gather*}
As is well-known, $\{m_{\lambda}(x)\}_{\lambda\in{\mathcal{P}}}$ form a~basis of $\Lambda_{N}(q,t)$.
Let $p_{n}(x):=\sum\limits_{i=1}^{N}x_{i}^{n}$, $n\in{\mathbf{Z}_{>0}}$ be the power sum, and for a partition~$\lambda$ we def\/ine $p_{\lambda}(x):=p_{\lambda_{1}}(x)\cdots p_{\lambda_{\ell(\lambda)}}(x)$.

Let $\rho^{N+1}_{N}: \Lambda_{N+1}(q,t) \to \Lambda_{N}(q,t)$ be the homomorphism def\/ined by
\begin{gather*}
\big(\rho^{N+1}_{N}f\big)(x_{1},\dots,x_{N}):=f(x_{1},\dots,x_{N},0),
\qquad
f\in\Lambda_{N+1}(q,t).
\end{gather*}
Def\/ine the ring of symmetric functions $\Lambda(q,t)$ as the projective limit def\/ined by $\big\{\rho^{N+1}_{N}\big\}_{N\geq{1}}$
\begin{gather*}
\Lambda(q,t):=\lim_{\longleftarrow}\Lambda_{N}(q,t).
\end{gather*}
It is known that $\{p_{\lambda}(x)\}_{\lambda\in{\mathcal{P}}}$ form a~basis of $\Lambda(q,t)$.
Then we def\/ine an inner product $\langle\;, \;\rangle_{q,t}$ as follows
\begin{gather*}
\langle p_{\lambda}(x), p_{\mu}(x) \rangle_{q,t}=\delta_{\lambda\mu}z_{\lambda}(q,t).
\end{gather*}

We def\/ine the order in $\mathcal{P}$ as follows: for $\lambda, \mu\in\mathcal{P}$,
\begin{gather*}
\lambda\geq{\mu}
\quad
\Longleftrightarrow
\quad
|\lambda|=|\mu|
\qquad
\text{and for all~}i
\qquad
\lambda_{1}+\dots+\lambda_{i}\geq{\mu_{1}+\dots+\mu_{i}}.
\end{gather*}
The existence of the Macdonald symmetric functions~\cite{4} is stated as follows: For each partition $\lambda$, there exists a unique symmetric function $P_{\lambda}(x)\in\Lambda(q,t)$ satisfying the following conditions:
\begin{gather*}
(1)
\quad
P_{\lambda}(x)=\sum\limits_{\mu\leq{\lambda}}u_{\lambda\mu}m_{\mu}(x),
\qquad
u_{\lambda\mu}\in{\mathbf{Q}(q,t)},
\\
(2)
\quad
\lambda \neq \mu \ \Longrightarrow \ \langle P_{\lambda}(x), P_{\mu}(x)\rangle_{q,t}=0.
\end{gather*}

For a~Macdonald symmetric function $P_{\lambda}(x)$, we def\/ine the $N$-variable symmetric polynomial
$P_{\lambda}(x_{1},\dots,x_{N})$ as $P_{\lambda}(x_{1},\dots,x_{N}):=P_{\lambda}(x_{1},\dots,x_{N},0,0,\dots)$,
$\ell(\lambda)\leq{N}$.
We call it the $N$-variables Macdonald polynomials.
We set the $q$-shift operator by
\begin{gather*}
T_{q,x_{i}}f(x_{1},\dots, x_{N}):=f(x_{1},\dots,qx_{i},\dots,x_{N})
\end{gather*}
and def\/ine the Macdonald operator $H_{N}(q,t): \Lambda_{N}(q,t) \to \Lambda_{N}(q,t)$ as follows
\begin{gather*}
H_{N}(q,t):=\sum\limits_{i=1}^{N}\prod\limits_{j\neq{i}}\frac{tx_{i}-x_{j}}{x_{i}-x_{j}}T_{q,x_{i}}.
\end{gather*}

Then for each partition $\lambda$, $\ell(\lambda)\leq{N}$,
the Macdonald polynomial $P_{\lambda}(x_{1},\dots,x_{N})$ is
an eigenfunction of the Macdonald operator~\cite{4}
\begin{gather*}
H_{N}(q,t)P_{\lambda}(x_{1},\dots,x_{N})=\varepsilon_{N}(\lambda)P_{\lambda}(x_{1},\dots,x_{N}),
\qquad
\varepsilon_{N}(\lambda):=\sum\limits_{i=1}^{N}q^{\lambda_{i}}t^{N-i}.
\end{gather*}

In the following $[f(z)]_{1}$ stands for the constant term of $f(z)$ in $z$.

\begin{Proposition}[free f\/ield realization of the Macdonald operator]\label{pr1.3}
Define the operator $\phi(z): \mathcal{F} \to \mathcal{F}\otimes\mathbf{C}[[z,z^{-1}]]$ as follows
\begin{gather*}
\phi(z):=\exp\left(\sum\limits_{n>0}\frac{1-t^{n}}{1-q^{n}}a_{-n}\frac{z^{n}}{n}\right),
\end{gather*}
and set $\phi_{N}(x):=\prod\limits_{j=1}^{N}\phi(x_{j})$.

$(1)$ The operator $\eta(z)$ reproduces the Macdonald operator $H_{N}(q,t)$ as follows
\begin{gather*}
[\eta(z)]_{1}\phi_{N}(x)|0 \rangle=t^{-N}\{(t-1)H_{N}(q,t)+1\}\phi_{N}(x)|0 \rangle.
\end{gather*}

$(2)$ The operator $\xi(z)$ reproduces the Macdonald operator $H_{N}(q^{-1},t^{-1})$ as follows
\begin{gather*}
[\xi(z)]_{1}\phi_{N}(x)|0 \rangle=t^{N}\big\{\big(t^{-1}-1\big)H_{N}\big(q^{-1},t^{-1}\big)+1\big\}\phi_{N}(x)|0 \rangle.
\end{gather*}
\end{Proposition}

We also have the dual version of Proposition~\ref{pr1.3}.

\begin{Proposition}[dual version of Proposition~\ref{pr1.3}]
Define the operator $\phi^{\ast}(z): \mathcal{F}^{\ast}\to\mathcal{F}^{\ast}\otimes\mathbf{C}[[z,z^{-1}]]$ by
\begin{gather*}
\phi^{\ast}(z):=\exp\left(\sum\limits_{n>0}\frac{1-t^{n}}{1-q^{n}}a_{n}\frac{z^{n}}{n}\right),
\end{gather*}
and set $\phi^{\ast}_{N}(x):=\prod\limits_{j=1}^{N}\phi^{\ast}(x_{j})$.

$(1)$ The operator $\eta(z)$ reproduces the Macdonald operator $H_{N}(q,t)$ as follows
\begin{gather*}
\langle 0|\phi^{\ast}_{N}(x)[\eta(z)]_{1}=t^{-N}\{(t-1)H_{N}(q,t)+1\}\langle 0|\phi^{\ast}_{N}(x).
\end{gather*}

$(2)$ The operator $\xi(z)$ reproduces the Macdonald operator $H_{N}(q^{-1},t^{-1})$ as follows
\begin{gather*}
\langle 0|\phi^{\ast}_{N}(x)[\xi(z)]_{1}=t^{N}\big\{\big(t^{-1}-1\big)H_{N}\big(q^{-1},t^{-1}\big)+1\big\}\langle 0|\phi^{\ast}_{N}(x).
\end{gather*}
\end{Proposition}

\begin{Remark}
Let us def\/ine the kernel function of the Macdonald operator
$\Pi_{MN}(q,t)(x,y)$, $M,N\in\mathbf{Z}_{>0}$,
  by
\begin{gather*}
\Pi_{MN}(q,t)(x,y)
:=\prod\limits_{\genfrac{}{}{0pt}{1}{1\leq{i}\leq{M}}{1\leq{j}\leq{N}}}\frac{(tx_{i}y_{j};q)_{\infty}}{(x_{i}y_{j};q)_{\infty}}.
\end{gather*}
Then the kernel function $\Pi_{MN}(q,t)(x,y)$ is reproduced from the operators
$\phi^{\ast}_{M}(x)$, $\phi_{N}(y)$ as
\begin{gather*}
\langle 0|\phi^{\ast}_{M}(x)\phi_{N}(y)|0 \rangle=\Pi_{MN}(q,t)(x,y).
\end{gather*}
\end{Remark}

\subsection[Trigonometric Feigin--Odesskii algebra $\mathcal{A}$]{Trigonometric Feigin--Odesskii algebra $\boldsymbol{\mathcal{A}}$}

In this subsection, we review basic facts of the
trigonometric Feigin--Odesskii algebra~\cite{1}.
In the following let $q,t\in\mathbf{C}$ be parameters satisfying $|q|<1$.

\begin{Definition}[trigonometric Feigin--Odesskii algebra $\mathcal{A}$]
\label{def1.6}
Let $\varepsilon_{n}(q;x)$, $n\in\mathbf{Z}_{>0}$ be a~function def\/ined as
\begin{gather*}
\varepsilon_{n}(q;x):=\prod\limits_{1\leq{a}<b\leq{n}}\frac{(x_{a}-qx_{b})(x_{a}-q^{-1}x_{b})}{(x_{a}-x_{b})^{2}}.
\end{gather*}
We also def\/ine $\omega(x,y)$ as
\begin{gather*}
\omega(x,y):=\frac{(x-q^{-1}y)(x-ty)(x-qt^{-1}y)}{(x-y)^{3}}.
\end{gather*}
For an $N$-variable function $f(x_{1},\dots,x_{N})$, we def\/ine the action of the symmetric group $\mathfrak{S}_{N}$ of order~$N$ on~$f(x_{1},\dots,x_{N})$ by
$\sigma\cdot(f(x_{1},\dots,x_{N})):=f(x_{\sigma(1)},\dots,x_{\sigma(N)})$, $\sigma\in\mathfrak{S}_{N}$.
We def\/ine the symmetrizer as
\begin{gather*}
\operatorname{Sym}[f(x_{1},\dots,x_{N})]:=\frac{1}{N!}\sum\limits_{\sigma\in\mathfrak{S}_{N}}\sigma\cdot(f(x_{1},\dots,x_{N})).
\end{gather*}
For an $m$-variable function $f(x_{1},\dots,x_{m})$ and an $n$-variable function $g(x_{1},\dots,x_{n})$, we def\/ine the
star product~$\ast$ as follows
\begin{gather*}
(f\ast g)(x_{1},\dots,x_{m+n}):=\operatorname{Sym}\left[f(x_{1},\dots,x_{m})g(x_{m+1},\dots,x_{m+n})
\prod\limits_{\genfrac{}{}{0pt}{1}{1\leq{\alpha}\leq{m}}{m+1\leq{\beta}\leq{m+n}}}\omega(x_{\alpha},x_{\beta})\right].
\end{gather*}
For a partition $\lambda$, we def\/ine $\varepsilon_{\lambda}(q;x)$, $x=(x_{1},\dots,x_{|\lambda|})$ as
\begin{gather*}
\varepsilon_{\lambda}(q;x):=(\varepsilon_{\lambda_{1}}(q;\bullet)\ast\cdots\ast\varepsilon_{\lambda_{\ell(\lambda)}}(q;\bullet))(x).
\end{gather*}

Set $\mathcal{A}_{0}:=\mathbf{Q}(q,t)$, $\mathcal{A}_{n}:=\operatorname{span}\{\varepsilon_{\lambda}(q;x): |\lambda|=n\}$,
$n\geq{1}$.
We def\/ine the trigonometric Feigin--Odesskii algebra to be $\mathcal{A}:=\bigoplus_{n\geq{0}}\mathcal{A}_{n}$ whose algebra structure is given by the star product~$\ast$.
\end{Definition}

\begin{Remark}
The def\/inition of the trigonometric Feigin--Odesskii algebra $\mathcal{A}$ above is a~short-handed one of~\cite{1}.
For instance, there would be a~question why the function $\varepsilon_{n}(q;x)$ appears.
For more detail of the trigonometric Feigin--Odesskii algebra $\mathcal{A}$, see~\cite{1}.
\end{Remark}

\begin{Proposition}[\cite{1}]\label{proposition2.8}
The trigonometric Feigin--Odesskii algebra $(\mathcal{A},\ast)$ is unital, associative, and commutative.
\end{Proposition}

\subsection[Commutative families $\mathcal{M}$, $\mathcal{M}^{\prime}$]{Commutative families $\boldsymbol{\mathcal{M}}$,
$\boldsymbol{\mathcal{M}^{\prime}}$}

Here we give an overview of the construction of
the commutative families of the Macdonald operator using the Ding--Iohara--Miki algebra and the trigonometric
Feigin--Odesskii algebra.

\begin{Definition}
[map $\mathcal{O}$] Def\/ine a linear map $\mathcal{O} : \mathcal{A} \to \text{End}(\mathcal{F})$ as
\begin{gather*}
\mathcal{O}(f):=\left[f(z_{1},\dots,z_{n})\prod_{1\leq{i<j}\leq{n}}\omega(z_{i},z_{j})^{-1}\eta(z_{1})\cdots\eta(z_{n})\right]_{1}
\end{gather*}
for $f\in\mathcal{A}_{n}$, where $[f(z_{1},\dots,z_{n})]_{1}$ denotes the constant term of $f(z_{1},\dots,z_{n})$ in $z_{1},\dots,z_{n}$, and extend to~$\mathcal{A}$ linearly.
\end{Definition}

From the relation
\begin{gather*}
\eta(z)\eta(w)=g\left(\frac{z}{w}\right)\eta(w)\eta(z),
\end{gather*}
we have the following
\begin{align*}
\frac{1}{\omega(z,w)}\eta(z)\eta(w)=\frac{1}{\omega(w,z)}\eta(w)\eta(z).
\end{align*}
This relation shows that the operator-valued function
\begin{gather*}
\prod\limits_{1\leq{i<j}\leq{N}}\omega(x_{i},x_{j})^{-1}\eta(x_{1})\cdots\eta(x_{N})
\end{gather*}
is symmetric in $x_{1},\dots,x_{N}$.
From this fact, we have the following proposition.

\begin{Proposition}
\label{pr1.10}
The map $\mathcal{O}$ and the star product~$\ast$ are compatible: for $f,g\in\mathcal{A}$, we have $\mathcal{O}(f\ast
g)=\mathcal{O}(f)\mathcal{O}(g)$.
\end{Proposition}

The trigonometric Feigin--Odesskii algebra $\mathcal{A}$ is commutative with respect to the star product~$\ast$, therefore
we have the following corollary.

Let $V$ be a~$\mathbf{C}$-vector space and $T: V \to V$ be a~$\mathbf{C}$-linear operator.
Then for a~subset $W\subset V$, the symbol $T|_{W}$ denotes the restriction of $T$ on $W$.
For a~subset $M\subset\operatorname{End}_{\mathbf{C}}(V)$, we use the symbol $M|_{W}:=\{T|_{W}: T\in M\}$.

\begin{Corollary}[commutative family $\mathcal{M}$]
\quad
\begin{enumerate}\itemsep=0pt
\item[$(1)$] Set $\mathcal{M}:=\mathcal{O}(\mathcal{A})$.
The space $\mathcal{M}$ consists of operators commuting:
$[\mathcal{O}(f),\mathcal{O}(g)]=0$, \mbox{$f,g\in\mathcal{A}$}.

\item[$(2)$] The space $\mathcal{M}|_{\mathbf{C}\phi_{N}(x)|0 \rangle}$ is a~set of commuting $q$-difference operators containing
the Macdonald operator $H_{N}(q,t)$.
\end{enumerate}
\end{Corollary}

\begin{proof}
(1) This statement follows from the commutativity of $\mathcal{A}$ in terms of the star product~$\ast$ and
Proposition~\ref{pr1.10}.

(2) Due to the free f\/ield realization of the Macdonald operator $H_{N}(q,t)$, the operator\linebreak
$\mathcal{O}(\varepsilon_{r}(q;z))$, $r\in\mathbf{Z}_{>0}$, acts on $\phi_{N}(x)|0 \rangle$, $N\in\mathbf{Z}_{>0}$,
as an $r$-th order $q$-dif\/ference operator.
By the fact that $\mathcal{M}=\mathcal{O}(\mathcal{A})$ is generated by
$\{\mathcal{O}(\varepsilon_{r}(q;z))\}_{r\in\mathbf{Z}_{>0}}$, and the relation
\begin{gather*}
\mathcal{O}(\varepsilon_{1}(q;z))\phi_{N}(x)|0 \rangle=[\eta(z)]_{1}\phi_{N}(x)|0
\rangle=t^{-N}\{(t-1)H_{N}(q,t)+1\}\phi_{N}(x)|0 \rangle
\end{gather*}
and (1), the statement holds.
\end{proof}

The Macdonald operator $H_{N}(q^{-1},t^{-1})$ is reproduced from the operator $\xi(z)$.
By this fact, we can construct another commutative family of the Macdonald operator.

\begin{Definition}[map $\mathcal{O}^{\prime}$]
Def\/ine a function $\omega^{\prime}(x,y)$ as
\begin{gather*}
\omega^{\prime}(x,y):=\frac{(x-qy)(x-t^{-1}y)(x-q^{-1}ty)}{(x-y)^{3}}.
\end{gather*}
Def\/ine a linear map $\mathcal{O}^{\prime} : \mathcal{A} \to \text{End}(\mathcal{F})$ as
\begin{gather*}
\mathcal{O}^{\prime}(f):=\left[f(z_{1},\dots,z_{n})\prod_{1\leq{i<j}\leq{n}}\omega^{\prime}(z_{i},z_{j})^{-1}\xi(z_{1})\cdots\xi(z_{n})\right]_{1}
\end{gather*}
for $f\in\mathcal{A}_{n}$, and extend to $\mathcal{A}$ linearly.
\end{Definition}

From the relation $\omega(x,y)=\omega^{\prime}(y,x)$ we have
\begin{Lemma}\label{lem1.13}
Define another star product $\ast^{\prime}$ as follows
\begin{gather*}
(f\ast^{\prime} g)(x_{1},\dots,x_{m+n}):=\operatorname{Sym}\left[f(x_{1},\dots,x_{m})g(x_{m+1},\dots,x_{m+n})
\prod\limits_{\genfrac{}{}{0pt}{1}{1\leq{\alpha}\leq{m}}{m+1\leq{\beta}\leq{m+n}}}\omega^{\prime}(x_{\alpha},x_{\beta})\right].
\end{gather*}
Then in the trigonometric Feigin--Odesskii algebra $\mathcal{A}$, we have $\ast^{\prime}=\ast$.
\end{Lemma}

We can check the map $\mathcal{O}^{\prime}$ and the star product $\ast^{\prime}$ are compatible in the similar way of
the proof of Proposition~\ref{pr1.10}.
Furthermore since $\ast^{\prime}=\ast$, we have

\begin{Corollary}[commutative family $\mathcal{M}^{\prime}$]\quad
\begin{enumerate}\itemsep=0pt
\item[$(1)$] Set $\mathcal{M}^{\prime}:=\mathcal{O}^{\prime}(\mathcal{A})$.
Then the space $\mathcal{M}^{\prime}$ consists of commuting operators.

\item[$(2)$] The space $\mathcal{M}^{\prime}|_{\mathbf{C}\phi_{N}(x)|0 \rangle}$ is a~set of commuting $q$-difference operators
containing the Macdonald operator $H_{N}(q^{-1},t^{-1})$.
\end{enumerate}
\end{Corollary}

From the relation $[[\eta(z)]_{1},[\xi(w)]_{1}]=0$, we have the following proposition.

\begin{Proposition}
\label{pr1.15}
The commutative families $\mathcal{M}$, $\mathcal{M}^{\prime}$ satisfy $[\mathcal{M},\mathcal{M}^{\prime}]=0$.
\end{Proposition}

\begin{proof}
This proposition follows from the existence of the Macdonald symmetric functions.
That is, elements of the commutative families are simultaneously diagonalized by the Macdonald symmetric functions.
\end{proof}

From Proposition~\ref{pr1.15}, commutative families $\mathcal{M}|_{\mathbf{C}\phi_{N}(x)|0 \rangle}$,
$\mathcal{M}^{\prime}|_{\mathbf{C}\phi_{N}(x)|0 \rangle}$ also commute:\linebreak
$[\mathcal{M}|_{\mathbf{C}\phi_{N}(x)|0 \rangle}, \mathcal{M}^{\prime}|_{\mathbf{C}\phi_{N}(x)|0 \rangle}]=0$.

\section{Elliptic case}\label{Sec2}

In this section, we construct commutative families of the elliptic Macdonald operators by using the elliptic Ding--Iohara--Miki algebra and the elliptic Feigin--Odesskii algebra.
Let $q, t, p\in\mathbf{C}$ with $|q|<1$, $|p|<1$.

\subsection[Elliptic Ding--Iohara--Miki algebra $\mathcal{U}(q,t,p)$]{Elliptic Ding--Iohara--Miki algebra $\boldsymbol{\mathcal{U}(q,t,p)}$}

The elliptic Ding--Iohara--Miki algebra is an
elliptic analog of the Ding--Iohara--Miki algebra introduced in~\cite{3}.
First we recall the def\/inition of the elliptic Ding--Iohara--Miki algebra and its free f\/ield realization.

\begin{Definition}
[elliptic Ding--Iohara--Miki algebra $\mathcal{U}(q,t,p)$] Def\/ine the structure function
$g_{p}(x)\in\mathbf{C}[[x,x^{-1}]]$ as
\begin{gather*}
g_{p}(x):=\frac{\Theta_{p}(qx)\Theta_{p}(t^{-1}x)\Theta_{p}(q^{-1}tx)}{\Theta_{p}(q^{-1}x)\Theta_{p}(tx)\Theta_{p}(qt^{-1}x)}.
\end{gather*}
Let $x^{\pm}(p;z):=\sum\limits_{n\in{\mathbf{Z}}}x^{\pm}_{n}(p)z^{-n}$,
$\psi^{\pm}(p;z):=\sum\limits_{n\in{\mathbf{Z}}}\psi^{\pm}_{n}(p)z^{-n}$ be currents and $C$ be a~central, inver\-tible
element satisfying the following relations
\begin{gather*}
[\psi^{\pm}(p;z), \psi^{\pm}(p;w)]=0,
\qquad
\psi^{+}(p;z)\psi^{-}(p;w)=\frac{g_{p}(Cz/w)}{g_{p}(C^{-1}z/w)}\psi^{-}(p;w)\psi^{+}(p;z),
\\
\psi^{\pm}(p;z)x^{+}(p;w)=g_{p}\left(C^{\pm\frac{1}{2}}\frac{z}{w}\right)x^{+}(p;w)\psi^{\pm}(p;z),
\\
\psi^{\pm}(p;z)x^{-}(p;w)=g_{p}\left(C^{\mp\frac{1}{2}}\frac{z}{w}\right)^{-1}x^{-}(p;w)\psi^{\pm}(p;z),
\\
x^{\pm}(p;z)x^{\pm}(p;w)=g_{p}\left(\frac{z}{w}\right)^{\pm 1}x^{\pm}(p;w)x^{\pm}(p;z),
\\
[x^{+}(p;z),x^{-}(p;w)]
=\frac{\Theta_{p}(q)\Theta_{p}(t^{-1})}{(p;p)_{\infty}^{3}\Theta_{p}(qt^{-1})}\\
\hphantom{[x^{+}(p;z),x^{-}(p;w)]=}{}\times
\left\{\delta\left(C\frac{w}{z}\right)\psi^{+}(p;C^{1/2}w)-\delta\left(C^{-1}\frac{w}{z}\right)\psi^{-}(p;C^{-1/2}w)\right\}.
\end{gather*}
We def\/ine the elliptic Ding--Iohara--Miki algebra $\mathcal{U}(q,t,p)$ to be the associative $\mathbf{C}$-algebra
gene\-ra\-ted by $\{x^{\pm}_{n}(p)\}_{n\in{\mathbf{Z}}}$, $\{\psi^{\pm}_{n}(p)\}_{n\in{\mathbf{Z}}}$ and $C$.
\end{Definition}

Let $\mathcal{B}_{a,b}$ be the associative $\mathbf{C}$-algebra generated by $\{a_{n}\}_{n\in{\mathbf{Z}\setminus\{0\}}}$, $\{b_{n}\}_{n\in{\mathbf{Z}\setminus\{0\}}}$ and the following relations
\begin{gather*}
[a_{m},a_{n}]=m(1-p^{|m|})\frac{1-q^{|m|}}{1-t^{|m|}}\delta_{m+n,0},
\qquad
[b_{m},b_{n}]=m\frac{1-p^{|m|}}{(qt^{-1}p)^{|m|}}\frac{1-q^{|m|}}{1-t^{|m|}}\delta_{m+n,0},
\\
[a_{m},b_{n}]=0, \qquad m,n\in\mathbf{Z}\setminus\{0\}.
\end{gather*}
We def\/ine the normal ordering $\bm{:} \bullet \bm{:}$ as usual
\begin{gather*}
\bm{:}a_{m}a_{n}\bm{:}=
\begin{cases}
a_{m}a_{n}, & m<n,
\\
a_{n}a_{m}, & m\geq{n},
\end{cases}
\qquad
\bm{:}b_{m}b_{n}\bm{:}=
\begin{cases}
b_{m}b_{n}, & m<n,
\\
b_{n}b_{m}, & m\geq{n}.
\end{cases}
\end{gather*}
Let $|0 \rangle$ be the vacuum vector which satisf\/ies the condition $a_{n}|0 \rangle=b_{n}|0 \rangle=0$, $n>0$, and set
the boson Fock space $\mathcal{F}$ as the left $\mathcal{B}_{a,b}$ module
\begin{gather*}
\mathcal{F}=\operatorname{span}\{a_{-\lambda}b_{-\mu}|0 \rangle: \lambda, \mu\in{\mathcal{P}}\}.
\end{gather*}
Let $\langle 0|$ be the dual vacuum vector which satisf\/ies the condition $\langle 0|a_{n}=\langle 0|b_{n}=0$, $n<0$, and
$\langle 0|a_{0}=0$.
We def\/ine the dual boson Fock space as the right $\mathcal{B}_{a,b}$ module
\begin{gather*}
\mathcal{F}^{\ast}:=\operatorname{span}\{\langle 0|a_{\lambda}b_{\mu}: \lambda,\mu\in\mathcal{P}\}.
\end{gather*}
For a~partition $\lambda$, set $n_{\lambda}(a)=\sharp\{i: \lambda_{i}=a\}$,
$z_{\lambda}=\prod\limits_{a\geq{1}}a^{n_{\lambda}(a)}n_{\lambda}(a)!$ and def\/ine $z_{\lambda}(q,t,p)$,
$\overline{z}_{\lambda}(q,t,p)$~by
\begin{gather*}
z_{\lambda}(q,t,p):=z_{\lambda}\prod\limits_{i=1}^{\ell(\lambda)}(1-p^{\lambda_{i}})\frac{1-q^{\lambda_{i}}}{1-t^{\lambda_{i}}},
\qquad
\overline{z}_{\lambda}(q,t,p):=z_{\lambda}\prod\limits_{i=1}^{\ell(\lambda)}
\frac{1-p^{\lambda_{i}}}{(qt^{-1}p)^{\lambda_{i}}}\frac{1-q^{\lambda_{i}}}{1-t^{\lambda_{i}}}.
\end{gather*}
We def\/ine a~bilinear form $\langle \bullet|\bullet \rangle: \mathcal{F}^{\ast} \times \mathcal{F} \to \mathbf{C}$ by the
following conditions:
\begin{gather*}
(1)
\quad
\langle 0|0 \rangle=1,
\\
(2)
\quad
\langle 0|a_{\lambda_{1}}b_{\lambda_{2}}a_{-\mu_{1}}b_{-\mu_{2}}|0
\rangle=\delta_{\lambda_{1}\mu_{1}}\delta_{\lambda_{2}\mu_{2}}z_{\lambda_{1}}(q,t,p)\overline{z}_{\lambda_{2}}(q,t,p).
\end{gather*}

\begin{Theorem}[free f\/ield realization of the elliptic Ding--Iohara--Miki algebra $\mathcal{U}(q,t,p)$]
\label{th2.2}
Define operators $\eta(p;z)$, $\xi(p;z)$, $\varphi^{\pm}(p;z): \mathcal{F} \to \mathcal{F}\otimes\mathbf{C}[[z,z^{-1}]]$
as follows $(\gamma:=(qt^{-1})^{-1/2})$
\begin{gather*}
\eta(p;z):=\bm{:}\exp\left(-\sum\limits_{n\neq{0}}\frac{1-t^{-n}}{1-p^{|n|}}p^{|n|}b_{n}\frac{z^{n}}{n}\right)
\exp\left(-\sum\limits_{n\neq{0}}\frac{1-t^{n}}{1-p^{|n|}}a_{n}\frac{z^{-n}}{n}\right)\bm{:},
\\
\xi(p;z):=\bm{:}\exp\left(\sum\limits_{n\neq{0}}\frac{1-t^{-n}}{1-p^{|n|}}\gamma^{-|n|}p^{|n|}b_{n}\frac{z^{n}}{n}\right)
\exp\left(\sum\limits_{n\neq{0}}\frac{1-t^{n}}{1-p^{|n|}}\gamma^{|n|}a_{n}\frac{z^{-n}}{n}\right)\bm{:},
\\
\varphi^{+}(p;z):=\bm{:}\eta(p;\gamma^{1/2}z)\xi(p;\gamma^{-1/2}z)\bm{:},
\qquad
\varphi^{-}(p;z):=\bm{:}\eta(p;\gamma^{-1/2}z)\xi(p;\gamma^{1/2}z)\bm{:} .
\end{gather*}
Then the map
\begin{gather*}
C \mapsto \gamma,
\qquad
x^{+}(p;z) \mapsto \eta(p;z),
\qquad
x^{-}(p;z) \mapsto \xi(p;z),
\qquad
\psi^{\pm}(p;z) \mapsto \varphi^{\pm}(p;z)
\end{gather*}
gives a~representation of the elliptic Ding--Iohara--Miki algebra $\mathcal{U}(q,t,p)$.
\end{Theorem}

The elliptic Macdonald operator $H_{N}(q,t,p)$, $N\in\mathbf{Z}_{>0}$, is def\/ined as follows
\begin{gather*}
H_{N}(q,t,p):=\sum\limits_{i=1}^{N}\prod\limits_{j\neq{i}}\frac{\Theta_{p}(tx_{i}/x_{j})}{\Theta_{p}(x_{i}/x_{j})}T_{q,x_{i}}.
\end{gather*}
By the operators $\eta(p;z)$, $\xi(p;z)$ in Theorem~\ref{th2.2}, we can reproduce the elliptic Macdonald ope\-ra\-tor as
follows~\cite{3}.

\begin{Theorem}[free f\/ield realization of the elliptic Macdonald operator]\label{th2.3}
Let us define an ope\-ra\-tor $\phi(p;z): \mathcal{F} \to \mathcal{F}\otimes\mathbf{C}[[z,z^{-1}]]$ as follows
\begin{gather*}
\phi(p;z):=\exp\left(\sum\limits_{n>0}\frac{(1-t^{n})(qt^{-1}p)^{n}}{(1-q^{n})(1-p^{n})}b_{-n}\frac{z^{-n}}{n}\right)
\exp\left(\sum\limits_{n>0}\frac{1-t^{n}}{(1-q^{n})(1-p^{n})}a_{-n}\frac{z^{n}}{n}\right),
\end{gather*}
and put $\phi_{N}(p;x):=\prod\limits_{j=1}^{N}\phi(p;x_{j})$.

$(1)$ The elliptic Macdonald operator $H_{N}(q,t,p)$ is reproduced by the operator $\eta(p;z)$ as follows
\begin{gather*}
[\eta(p;z)-t^{-N}(\eta(p;z))_{-}(\eta(p;p^{-1}z))_{+}]_{1}\phi_{N}(p;x)|0 \rangle
=\frac{t^{-N+1}\Theta_{p}(t^{-1})}{(p;p)_{\infty}^{3}}H_{N}(q,t,p)\phi_{N}(p;x)|0 \rangle.
\end{gather*}
Here we use the notation $(\eta(p;z))_{\pm}$ as
\begin{gather*}
(\eta(p;z))_{\pm}:=\exp\left(-\sum\limits_{\pm
n>0}\frac{1-t^{-n}}{1-p^{|n|}}p^{|n|}b_{n}\frac{z^{n}}{n}\right)\exp\left(-\sum\limits_{\pm
n>0}\frac{1-t^{n}}{1-p^{|n|}}a_{n}\frac{z^{-n}}{n}\right).
\end{gather*}

$(2)$ The elliptic Macdonald operator $H_{N}(q^{-1},t^{-1},p)$ is reproduced by the operator $\xi(p;z)$ as follows
\begin{gather*}
[\xi(p;z)-t^{N}(\xi(p;z))_{-}(\xi(p;p^{-1}z))_{+}]_{1}\phi_{N}(p;x)|0 \rangle
=\frac{t^{N-1}\Theta_{p}(t)}{(p;p)_{\infty}^{3}}H_{N}\big(q^{-1},t^{-1},p\big)\phi_{N}(p;x)|0 \rangle.
\end{gather*}
Here we use the notation $(\xi(p;z))_{\pm}$ as
\begin{gather*}
(\xi(p;z))_{\pm}:=\exp\left(\sum\limits_{\pm
n>0}\frac{1-t^{-n}}{1-p^{|n|}}\gamma^{-|n|}p^{|n|}b_{n}\frac{z^{n}}{n}\right)\exp\left(\sum\limits_{\pm
n>0}\frac{1-t^{n}}{1-p^{|n|}}\gamma^{|n|}a_{n}\frac{z^{-n}}{n}\right).
\end{gather*}
\end{Theorem}

To state the next theorem, we introduce zero mode generators $a_{0}$, $Q$ satisfying
\begin{gather*}
[a_{0},Q]=1,
\qquad
[a_{n},a_{0}]=[b_{n},a_{0}]=0,
\qquad
[a_{n},Q]=[b_{n},Q]=0,
\qquad
n\in\mathbf{Z}\setminus\{0\}.
\end{gather*}
We also set the condition $a_{0}|0 \rangle=0$.
For a~complex number $\alpha\in\mathbf{C}$, we def\/ine $|\alpha \rangle:=e^{\alpha Q}|0 \rangle$.
Then we can check $a_{0}|\alpha \rangle=\alpha|\alpha \rangle$.
For $\alpha\in\mathbf{C}$, we set $\mathcal{F}_{\alpha}:=\operatorname{span}\{a_{-\lambda}b_{-\mu}|\alpha \rangle: \lambda,
\mu\in\mathcal{P}\}$.

\begin{Theorem}\label{th2.4}
Set $\widetilde{\eta}(p;z):=(\eta(p;z))_{-}(\eta(p;p^{-1}z))_{+}$,
$\widetilde{\xi}(p;z):=(\xi(p;z))_{-}(\xi(p;p^{-1}z))_{+}$ and define
\begin{gather*}
E(p;z):=\eta(p;z)-\widetilde{\eta}(p;z)t^{-a_{0}},
\qquad
F(p;z):=\xi(p;z)-\widetilde{\xi}(p;z)t^{a_{0}}.
\end{gather*}
Then the elliptic Macdonald operators $H_{N}(q,t,p)$, $H_{N}(q^{-1},t^{-1},p)$ are reproduced by the operators $E(p;z)$,
$F(p;z)$ as follows
\begin{gather*}
[E(p;z)]_{1}\phi_{N}(p;x)|N\rangle=\frac{t^{-N+1}\Theta_{p}(t^{-1})}{(p;p)_{\infty}^{3}}H_{N}(q,t,p)\phi_{N}(p;x)|N\rangle,
\\
[F(p;z)]_{1}\phi_{N}(p;x)|N\rangle=\frac{t^{N-1}\Theta_{p}(t)}{(p;p)_{\infty}^{3}}H_{N}\big(q^{-1},t^{-1},p\big)\phi_{N}(p;x)|N\rangle.
\end{gather*}
\end{Theorem}

Dual versions of Theorems~\ref{th2.3} and~\ref{th2.4} are also available.
For $\alpha\in\mathbf{C}$, set $\langle \alpha|:=\langle 0|e^{-\alpha Q}$.
Then we have $\langle \alpha|a_{0}=\alpha\langle \alpha|$.
For $\alpha\in\mathbf{C}$, we set $\mathcal{F}^{\ast}_{\alpha}:=\operatorname{span}\{\langle \alpha|a_{\lambda}b_{\mu}: \lambda,
\mu\in\mathcal{P}\}$.

\begin{Theorem}[dual versions of Theorems~\ref{th2.3} and~\ref{th2.4}]\label{th2.5}
Let us define the operator $\phi^{\ast}(p;z) : \mathcal{F}^{\ast} \to \mathcal{F}^{\ast}\otimes\mathbf{C}[[z,z^{-1}]]$ as follows
\begin{gather*}
\phi^{\ast}(p;z):=\exp\left(\sum\limits_{n>0}\frac{(1-t^{n})(qt^{-1}p)^{n}}{(1-q^{n})(1-p^{n})}b_{n}\frac{z^{-n}}{n}\right)
\exp\left(\sum\limits_{n>0}\frac{1-t^{n}}{(1-q^{n})(1-p^{n})}a_{n}\frac{z^{n}}{n}\right),
\end{gather*}
and set $\phi^{\ast}_{N}(p;x):=\prod\limits_{j=1}^{N}\phi^{\ast}(p;x_{j})$.

$(1)$ The elliptic Macdonald operators $H_{N}(q,t,p)$, $H_{N}(q^{-1},t^{-1},p)$ are reproduced by the ope\-ra\-tors
$\eta(p;z)$, $\xi(p;z)$ as follows
\begin{gather*}
\langle 0|\phi^{\ast}_{N}(p;x)\big[\eta(p;z)-t^{-N}(\eta(p;z))_{-}(\eta(p;p^{-1}z))_{+}\big]_{1}
=\frac{t^{-N+1}\Theta_{p}(t^{-1})}{(p;p)_{\infty}^{3}}H_{N}(q,t,p)\langle 0|\phi^{\ast}_{N}(p;x),
\\
\langle 0|\phi^{\ast}_{N}(p;x)\big[\xi(p;z)-t^{N}(\xi(p;z))_{-}(\xi(p;p^{-1}z))_{+}\big]_{1}
=\frac{t^{N-1}\Theta_{p}(t)}{(p;p)_{\infty}^{3}}H_{N}\big(q^{-1},t^{-1},p\big)\langle 0|\phi^{\ast}_{N}(p;x).
\end{gather*}

$(2)$ The operators $E(p;z)$, $F(p;z)$ reproduce the elliptic Macdonald operators $H_{N}(q,t,p)$, $H_{N}(q^{-1},t^{-1},p)$
as follows
\begin{gather*}
\langle N|\phi^{\ast}_{N}(p;x)[E(p;z)]_{1}=\frac{t^{-N+1}\Theta_{p}(t^{-1})}{(p;p)_{\infty}^{3}}H_{N}(q,t,p)
\langle N|\phi^{\ast}_{N}(p;x),
\\
\langle N|\phi^{\ast}_{N}(p;x)[F(p;z)]_{1}=\frac{t^{N-1}\Theta_{p}(t)}{(p;p)_{\infty}^{3}}H_{N}\big(q^{-1},t^{-1},p\big)
\langle N|\phi^{\ast}_{N}(p;x).
\end{gather*}
\end{Theorem}

\begin{Remark}
Let $\Pi_{MN}(q,t,p)(x,y)$, $M,N\in\mathbf{Z}_{>0}$, be the kernel function of the
elliptic Macdonald operator def\/ined as
\begin{gather*}
\Pi_{MN}(q,t,p)(x,y):=\prod\limits_{\genfrac{}{}{0pt}{1}{1\leq{i}\leq{M}}{1\leq{j}\leq{N}}}
\frac{\Gamma_{q,p}(x_{i}y_{j})}{\Gamma_{q,p}(tx_{i}y_{j})}.
\end{gather*}
Then the kernel function $\Pi_{MN}(q,t,p)(x,y)$ is reproduced from the operators
$\phi^{\ast}_{M}(p;x)$, $\phi_{N}(p;y)$ as
\begin{gather*}
\langle 0|\phi^{\ast}_{M}(p;x)\phi_{N}(p;y)|0 \rangle=\Pi_{MN}(q,t,p)(x,y).
\end{gather*}
\end{Remark}

\subsection[Elliptic Feigin--Odesskii algebra $\mathcal{A}(p)$]{Elliptic Feigin--Odesskii algebra $\boldsymbol{\mathcal{A}(p)}$}

The elliptic Feigin$-$Odesskii algebra is def\/ined quite similar as in the trigonometric case except for the emergence of elliptic functions~\cite{1}.
Let $q, t, p\in\mathbf{C}$ be complex parameters satifying $|q|<1$, $|p|<1$.

\begin{Definition}
[elliptic Feigin--Odesskii algebra $\mathcal{A}(p)$]
\label{def2.7}
Def\/ine an $n$-variable function $\varepsilon_{n}(q,p;x)$, $n\in\mathbf{Z}_{>0}$, as follows
\begin{gather*}
\varepsilon_{n}(q,p;x):=\prod\limits_{1\leq{a}<b\leq{n}}
\frac{\Theta_{p}(qx_{a}/x_{b})\Theta_{p}(q^{-1}x_{a}/x_{b})}{\Theta_{p}(x_{a}/x_{b})^{2}}.
\end{gather*}
Def\/ine a~function $\omega_{p}(x,y)$ as
\begin{gather*}
\omega_{p}(x,y):=\frac{\Theta_{p}(q^{-1}y/x)\Theta_{p}(ty/x)\Theta_{p}(qt^{-1}y/x)}{\Theta_{p}(y/x)^{3}}.
\end{gather*}
Def\/ine the star product~$\ast$ as
\begin{gather*}
(f\ast g)(x_{1},\dots,x_{m+n}):=\operatorname{Sym}\left[f(x_{1},\dots,x_{m})g(x_{m+1},\dots,x_{m+n})
\prod\limits_{\genfrac{}{}{0pt}{1}{1\leq{\alpha}\leq{m}}{m+1\leq{\beta}\leq{m+n}}}\omega_{p}(x_{\alpha},x_{\beta})\right].
\end{gather*}
For a partition $\lambda$, we set $\varepsilon_{\lambda}(q,p;x)$, $x=(x_{1},\dots,x_{|\lambda|})$ as
\begin{gather*}
\varepsilon_{\lambda}(q,p;x):=(\varepsilon_{\lambda_{1}}(q,p;\bullet)\ast\cdots\ast\varepsilon_{\lambda_{\ell(\lambda)}}(q,p;\bullet))(x).
\end{gather*}
Set $\mathcal{A}_{0}(p){:}=\mathbf{C}$, $\mathcal{A}_{n}(p){:}=\operatorname{span}\{\varepsilon_{\lambda}(q,p;x):
|\lambda|{=}n\}$, $n{\geq}{1}$.
We def\/ine the elliptic Feigin--Odesskii algebra as $\mathcal{A}(p):=\bigoplus_{n\geq{0}}\mathcal{A}_{n}(p)$ whose
algebra  structure is given by the star product~$\ast$.
\end{Definition}

As in the trigonometric case (Proposition~\ref{proposition2.8}), we have

\begin{Proposition}
The elliptic Feigin--Odesskii algebra $(\mathcal{A}(p),\ast)$ is an unital, associative, commutative algebra.
\end{Proposition}

\subsection[Commutative families $\mathcal{M}(p)$, $\mathcal{M}^{\prime}(p)$]{Commutative families
$\boldsymbol{\mathcal{M}(p)}$, $\boldsymbol{\mathcal{M}^{\prime}(p)}$}

For the operators $E(p;z)$, $F(p;z)$ in
Theorem~\ref{th2.4}, we have~\cite{3}

\begin{Proposition}\label{pr2.9}
The following relations hold
\begin{gather}
E(p;z)E(p;w)=g_{p}\left(\frac{z}{w}\right)E(p;w)E(p;z), \label{eq2.19}\\
F(p;z)F(p;w)=g_{p}\left(\frac{z}{w}\right)^{-1}F(p;w)F(p;z), \label{eq2.20}\\
[E(p;z),F(p;w)]=\frac{\Theta_{p}(q)\Theta_{p}(t^{-1})}{(p;p)_{\infty}^{3}
\Theta_{p}(qt^{-1})}\delta\left(\gamma\frac{w}{z}\right)\big\{\varphi^{+}\big(p;\gamma^{1/2}w\big)-\varphi^{+}
\big(p;\gamma^{1/2}p^{-1}w\big)\big\}.\label{eq2.21}
\end{gather}
\end{Proposition}

From the relation~\eqref{eq2.21} we have $[[E(p;z)]_{1},[F(p;w)]_{1}]=0$.
This corresponds to the commutativity of the elliptic Macdonald operators $[H_{N}(q,t,p),H_{N}(q^{-1},t^{-1},p)]=0$.

Def\/ine a~function $\omega_{p}^{\prime}(x,y)$ as
\begin{gather*}
\omega_{p}^{\prime}(x,y):=\frac{\Theta_{p}(qy/x)\Theta_{p}(t^{-1}y/x)\Theta_{p}(q^{-1}ty/x)}{\Theta_{p}(y/x)^{3}}.
\end{gather*}
Due to the relations~\eqref{eq2.19} and \eqref{eq2.20}, operator-valued functions
\begin{gather*}
\prod\limits_{1\leq{i<j}\leq{N}}\omega_{p}(x_{i},x_{j})^{-1}E(p;x_{1})\cdots E(p;x_{N}),
\qquad
\prod\limits_{1\leq{i<j}\leq{N}}\omega^{\prime}_{p}(x_{i},x_{j})^{-1}F(p;x_{1})\cdots F(p;x_{N})
\end{gather*}
are symmetric in $x_{1},\dots,x_{N}$.

\begin{Definition}[map $\mathcal{O}_{p}$] We def\/ine a linear map $\mathcal{O}_{p} : \mathcal{A}(p) \to \text{End}(\mathcal{F}_{\alpha})$, $\alpha\in\mathbf{C}$ as follows
\begin{gather*}
\mathcal{O}_{p}(f):=\left[f(z_{1},\dots,z_{n})\prod_{1\leq{i<j}\leq{n}}\omega_{p}(z_{i},z_{j})^{-1}E(p;z_{1})\cdots E(p;z_{n})\right]_{1}
\end{gather*}
for $f{\in}\mathcal{A}_{n}(p)$, where $[f(z_{1},\dots,z_{n})]_{1}$ denotes the constant term of $f(z_{1},\dots,z_{n})$ in $z_{1},\dots,z_{n}$, and extend linearly to~$\mathcal{A}(p)$.
\end{Definition}

In the similar way of the trigonometric case, we can check the following

\begin{Proposition}
The map $\mathcal{O}_{p}$ and the star product~$\ast$ are compatible: for $f,g\in\mathcal{A}(p)$, we have
$\mathcal{O}_{p}(f\ast g)=\mathcal{O}_{p}(f)\mathcal{O}_{p}(g)$.
\end{Proposition}

\begin{Theorem}[commutative family $\mathcal{M}(p)$]  \quad
\begin{enumerate}\itemsep=0pt
\item[$(1)$] Set $\mathcal{M}(p):=\mathcal{O}_{p}(\mathcal{A}(p))$.
The space $\mathcal{M}(p)$ is commutative.

\item[$(2)$] The space $\mathcal{M}(p)|_{\mathbf{C}\phi_{N}(p;x)|N \rangle}$ is a~set of commuting elliptic $q$-difference
operators containing the elliptic Macdonald operator $H_{N}(q,t,p)$.
\end{enumerate}
\end{Theorem}

A commutative family containing the elliptic Macdonald operator $H_{N}(q^{-1},t^{-1},p)$ is also constructed as follows

\begin{Definition}[map $\mathcal{O}_{p}^{\prime}$]
We def\/ine a linear map $\mathcal{O}_{p}^{\prime} : \mathcal{A}(p) \to \text{End}(\mathcal{F}_{\alpha})$, $\alpha\in\mathbf{C}$ as follows.
\begin{gather*}
\mathcal{O}_{p}^{\prime}(f):=\left[f(z_{1},\dots,z_{n})\prod_{1\leq{i<j}\leq{n}}\omega_{p}^{\prime}(z_{i},z_{j})^{-1}F(p;z_{1})\cdots F(p;z_{n})\right]_{1}
\end{gather*}
for $f\in\mathcal{A}_{n}(p)$, and extend linearly to~$\mathcal{A}(p)$.
\end{Definition}

As in the trigonometric case we have

\begin{Lemma}
Define another star product $\ast^{\prime}$ as
\begin{gather*}
(f\ast^{\prime} g)(x_{1},\dots,x_{m+n}):=\operatorname{Sym}\left[f(x_{1},\dots,x_{m})g(x_{m+1},\dots,x_{m+n})
\!\prod\limits_{\genfrac{}{}{0pt}{1}{1\leq{\alpha}\leq{m}}{m+1\leq{\beta}\leq{m+n}}}\!
\omega^{\prime}_{p}(x_{\alpha},x_{\beta})\right] .
\end{gather*}
In the elliptic Feigin--Odesskii algebra $\mathcal{A}(p)$, we have $\ast^{\prime}=\ast$.
\end{Lemma}

{\samepage \begin{Theorem}
[commutative family $\mathcal{M}^{\prime}(p)$]\quad
\begin{enumerate}\itemsep=0pt
\item[$(1)$] Set $\mathcal{M}^{\prime}(p):=\mathcal{O}_{p}^{\prime}(\mathcal{A}(p))$.
The space $\mathcal{M}^{\prime}(p)$ is commutative.

\item[$(2)$] The space $\mathcal{M}^{\prime}(p)|_{\mathbf{C}\phi_{N}(p;x)|N \rangle}$ is a~set of commuting elliptic
$q$-difference operators containing the elliptic Macdonald operator $H_{N}(q^{-1},t^{-1},p)$.
\end{enumerate}
\end{Theorem}}

Similar to Proposition~\ref{pr1.15}, we can show that the commutative families $\mathcal{M}(p)$,
$\mathcal{M}^{\prime}(p)$ commute with each other.

\begin{Theorem}
\label{th2.16}
Commutative families $\mathcal{M}(p)$, $\mathcal{M}^{\prime}(p)$ commute:
$[\mathcal{M}(p),\mathcal{M}^{\prime}(p)]{=}0$.
\end{Theorem}

Theorem~\ref{th2.16} is an elliptic analog of Proposition~\ref{pr1.15}.
But we can't prove Theorem~\ref{th2.16} in the similar way of the proof of Proposition~\ref{pr1.15}, because we don't
have an elliptic analog of the Macdonald symmetric functions.
Hence we will show Theorem~\ref{th2.16} in a~direct way.
We need the following lemma

\begin{Lemma}
\label{lem2.17}
Assume that an $r$-variable function $A(x_{1},\dots,x_{r})$ and an $s$-variable function $B(x_{1},\dots,x_{s})$ have
a~period $p$, i.e.
\begin{gather*}
T_{p,x_{i}}A(x_{1},\dots,x_{r})=A(x_{1},\dots,x_{r}),
\qquad
1\leq{i}\leq{r},
\\
T_{p,x_{i}}B(x_{1},\dots,x_{s})=B(x_{1},\dots,x_{s}),
\qquad
1\leq{i}\leq{s}.
\end{gather*}
Then we have
\begin{gather*}
\big[[A(z_{1},\dots,z_{r})E(p;z_{1})\cdots E(p;z_{r})]_{1},[B(w_{1},\dots,w_{s})F(p;w_{1})\cdots F(p;w_{s})]_{1}\big]=0.
\end{gather*}
\end{Lemma}

\begin{proof}
Recall the general formula of commutator
\begin{gather}
[A_{1}\cdots A_{r},B_{1}\cdots B_{s}]=\sum\limits_{i=1}^{r}\sum\limits_{j=1}^{s}A_{1}\cdots A_{i-1}B_{1}\cdots
B_{j-1}[A_{i},B_{j}]B_{j+1}\cdots B_{s}A_{i+1}\cdots A_{r}.
\label{eq2.27}
\end{gather}
Set $c(q,t,p):=\Theta_{p}(q)\Theta_{p}(t^{-1})/(p;p)_{\infty}^{3}\Theta_{p}(qt^{-1})$ and let $\Delta_{p}f(z):=f(pz)-f(z)$ stands for the $p$-dif\/ference of~$f(z)$.
By the identities~\eqref{eq2.21} and \eqref{eq2.27}, we have the following
\begin{gather}
[A(z_{1},\dots,z_{r})E(p;z_{1})\cdots E(p;z_{r}),B(w_{1},\dots,w_{s})F(p;w_{1})\cdots F(p;w_{s})]
\nonumber
\\
\qquad
=\sum\limits_{i=1}^{r}\sum\limits_{j=1}^{s}A(z_{1},\dots,z_{r})B(w_{1},\dots,w_{s})E(p;z_{1})\cdots
E(p;z_{i-1})F(p;w_{1})\cdots F(p;w_{j-1})
\nonumber
\\
\qquad
\phantom{=}
\times [E(p;z_{i}),F(p;w_{j})]F(p;w_{j+1})\cdots F(p;w_{s})E(p;z_{i+1})\cdots E(p;z_{r})
\nonumber
\\
\qquad
=c(q,t,p)\sum\limits_{i=1}^{r}\sum\limits_{j=1}^{s}E(p;z_{1})\cdots E(p;z_{i-1})F(p;w_{1})\cdots F(p;w_{j-1})
\nonumber
\\
\qquad
\phantom{=}
\times A(z_{1},\dots,\overbrace{\gamma w_{j}}^{i\text{-th}},\dots,z_{r})
B(w_{1},\dots,\overbrace{w_{j}}^{j\text{-th}},\dots,w_{s})
\delta\left(\gamma\frac{w_{j}}{z_{i}}\right)\Delta_{p}\varphi^{+}(p;\gamma^{1/2}p^{-1}w_{j})
\nonumber
\\
\qquad
\phantom{=}
\times F(p;w_{j+1})\cdots F(p;w_{s})E(p;z_{i+1})\cdots E(p;z_{r}).
\label{eq2.28}
\end{gather}
By picking up the constant term of $z_{i}$, $w_{j}$ dependent part of~\eqref{eq2.28}, we have
\begin{gather*}
\bigg[A(z_{1},\dots,\overbrace{\gamma w_{j}}^{i\text{-th}},\dots,z_{r})
B(w_{1},\dots,\overbrace{w_{j}}^{j\text{-th}},\dots,w_{s})
\delta\left(\gamma\frac{w_{j}}{z_{i}}\right)\Delta_{p}\varphi^{+}(p;\gamma^{1/2}p^{-1}w_{j})\bigg]_{z_{i},w_{j},1}
\\
\qquad
=\bigg[A(z_{1},\dots,\overbrace{\gamma w_{j}}^{i\text{-th}},\dots,z_{r})
B(w_{1},\dots,\overbrace{w_{j}}^{j\text{-th}},\dots,w_{s})\Delta_{p}\varphi^{+}(p;\gamma^{1/2}p^{-1}w_{j})\bigg]_{w_{j},1}
\\
\qquad
=\bigg[A(z_{1},\dots,\overbrace{\gamma w_{j}}^{i\text{-th}},\dots,z_{r})
B(w_{1},\dots,\overbrace{w_{j}}^{j\text{-th}},\dots,w_{s})\varphi^{+}(p;\gamma^{1/2}w_{j})\bigg]_{w_{j},1}
\\
\qquad
\phantom{=}
{}-\bigg[A(z_{1},\dots,\overbrace{\gamma w_{j}}^{i\text{-th}},\dots,z_{r})
B(w_{1},\dots,\overbrace{w_{j}}^{j\text{-th}},\dots,w_{s})\varphi^{+}(p;\gamma^{1/2}p^{-1}w_{j})\bigg]_{w_{j},1}.
\end{gather*}
We recall $[f(z)]_{1}\!=\![f(az)]_{1}$, $a\!\in\!\mathbf{C}$, and the assumption that both  $A(z_{1},\dots,z_{r})$ and $B(w_{1},\dots,w_{s})$ have the period~$p$. Therefore we have
\begin{gather*}
\bigg[A(z_{1},\dots,\overbrace{\gamma
w_{j}}^{i\text{-th}},\dots,z_{r})B(w_{1},\dots,\overbrace{w_{j}}^{j\text{-th}},\dots,w_{s})
\delta\left(\gamma\frac{w_{j}}{z_{i}}\right)\Delta_{p}\varphi^{+}(p;\gamma^{1/2}p^{-1}w_{j})\bigg]_{z_{i},w_{j},1}=0
\end{gather*}
for any~$i$,~$j$, proving the lemma.
\end{proof}

\begin{proof}[Proof of Theorem~\ref{th2.16}.] It is enough to show
\begin{gather*}
[\mathcal{O}_{p}(\varepsilon_{r}(q,p;z)),\mathcal{O}^{\prime}_{p}(\varepsilon_{s}(q,p;w))]=0,
\qquad
r,s\in\mathbf{Z}_{>0}.
\end{gather*}
By the def\/inition of $\mathcal{O}_{p}$, $\mathcal{O}^{\prime}_{p}$, the operators $\mathcal{O}_{p}(\varepsilon_{r}(q,p;z))$,
$\mathcal{O}^{\prime}_{p}(\varepsilon_{s}(q,p;w))$ are the constant terms of the following operators
\begin{gather}
\varepsilon_{r}(q,p;z)\prod\limits_{1\leq{i<j}\leq{r}}\omega_{p}(z_{i},z_{j})^{-1}E(p;z_{1})\cdots E(p;z_{r}),
\label{eq2.30}
\\
\varepsilon_{s}(q,p;w)\prod\limits_{1\leq{i<j}\leq{s}}\omega^{\prime}_{p}(w_{i},w_{j})^{-1}F(p;w_{1})\cdots F(p;w_{s}).
\label{eq2.31}
\end{gather}
Then their functional parts take the following forms
\begin{gather}
\text{(Functional part of~\eqref{eq2.30})}
\nonumber
\\
\qquad
=\varepsilon_{r}(q,p;z)\prod\limits_{1\leq{i<j}\leq{r}}\omega_{p}(z_{i},z_{j})^{-1}
=\prod\limits_{1\leq{i<j}\leq{r}}
\frac{\Theta_{p}(z_{i}/z_{j})\Theta_{p}(q^{-1}z_{i}/z_{j})}{\Theta_{p}(t^{-1}z_{i}/z_{j})\Theta_{p}(q^{-1}tz_{i}/z_{j})},
\label{eq2.32}
\\
\text{(Functional part of~\eqref{eq2.31})}
\nonumber
\\
\qquad
=\varepsilon_{s}(q,p;w)\prod\limits_{1\leq{i<j}\leq{s}}\omega^{\prime}_{p}(w_{i},w_{j})^{-1}
=\prod\limits_{1\leq{i<j}\leq{s}}
\frac{\Theta_{p}(w_{i}/w_{j})\Theta_{p}(qw_{i}/w_{j})}{\Theta_{p}(tw_{i}/w_{j})\Theta_{p}(qt^{-1}w_{i}/w_{j})}.
\label{eq2.33}
\end{gather}
We can check~\eqref{eq2.32},~\eqref{eq2.33} have a~period $p$.
By Lemma~\ref{lem2.17}, we have Theorem~\ref{th2.16}.
\end{proof}

By Theorem~\ref{th2.16}, commutative families $\mathcal{M}(p)|_{\mathbf{C}\phi_{N}(p;x)|N \rangle}$,
$\mathcal{M}^{\prime}(p)|_{\mathbf{C}\phi_{N}(p;x)|N \rangle}$ also commute:
$[\mathcal{M}(p)|_{\mathbf{C}\phi_{N}(p;x)|N \rangle}, \mathcal{M}^{\prime}(p)|_{\mathbf{C}\phi_{N}(p;x)|N \rangle}]=0$.

In the trigonometric case, relations between the commutative families $\mathcal{M}$, $\mathcal{M}^{\prime}$ and the higher-order Macdonald operators are studied in~\cite{1}. In contrast, the elliptic case relations between the commutative families $\mathcal{M}(p)$, $\mathcal{M}^{\prime}(p)$ and the higher-order elliptic Macdonald operator~\cite{2} still remain unclear.

\appendix

\section{Appendix}

\subsection{Trigonometric kernel function and its functional equation}

The following theorem is shown by Komori, Noumi, and Shiraishi~\cite{2}.

\begin{Theorem}
[\protect{\cite{2}}]
\label{theorem3.1}
Define the Macdonald operator $H_{N}(q,t)$, $N\in\mathbf{Z}_{>0}$, as
\begin{gather*}
H_{N}(q,t):=\sum\limits_{i=1}^{N}\prod\limits_{j\neq{i}}\frac{tx_{i}-x_{j}}{x_{i}-x_{j}}T_{q,x_{i}},
\qquad
T_{q,x}f(x):=f(qx),
\end{gather*}
and its kernel function $\Pi_{MN}(q,t)(x,y)$, $M,N\in\mathbf{Z}_{>0}$, as
\begin{gather*}
\Pi_{MN}(q,t)(x,y)
:=\prod\limits_{\genfrac{}{}{0pt}{1}{1\leq{i}\leq{M}}{1\leq{j}\leq{N}}}\frac{(tx_{i}y_{j};q)_{\infty}}{(x_{i}y_{j};q)_{\infty}}.
\end{gather*}
Then we have the following functional equation
\begin{gather}
\big\{H_{M}(q,t)_{x}{-}t^{M-N}H_{N}(q,t)_{y}\big\}\Pi_{MN}(q,t)(x,y)
=\frac{1{-}t^{M-N}}{1{-}t}\Pi_{MN}(q,t)(x,y).
\label{eq3.3}
\end{gather}
Here we denote the Macdonald operator which acts on functions of $x_{1},\dots,x_{M}$ by $H_{M}(q,t)_{x}$.
\end{Theorem}

In the following, we will show an elliptic analog of Theorem~\ref{theorem3.1} by the free f\/ield realization of the
elliptic Macdonald operator.

\subsection{Elliptic kernel function and its functional equation}

By the free f\/ield realization, we can show the following theorem.

\begin{Theorem}[functional equation of the elliptic kernel function]\label{theorem3.2}
Define the elliptic kernel function $\Pi_{MN}(q,t,p)(x,y)$ by
\begin{gather*}
\Pi_{MN}(q,t,p)(x,y)
:=\prod\limits_{\genfrac{}{}{0pt}{1}{1\leq{i}\leq{M}}{1\leq{j}\leq{N}}}\frac{\Gamma_{q,p}(x_{i}y_{j})}{\Gamma_{q,p}(tx_{i}y_{j})}.
\end{gather*}
We also define $C_{MN}(p;x,y)$ as
\begin{gather*}
C_{MN}(p;x,y):=\frac{\langle 0|\phi^{\ast}_{M}(p;x)[(\eta(p;z))_{-}(\eta(p;p^{-1}z))_{+}]_{1}\phi_{N}(p;y)|0
\rangle}{\Pi_{MN}(q,t,p)(x,y)}
\\
\hphantom{C_{MN}(p;x,y)}{}
=\left[\prod\limits_{i=1}^{M}\frac{\Theta_{p}(t^{-1}x_{i}z)}{\Theta_{p}(x_{i}z)}
\prod\limits_{j=1}^{N}\frac{\Theta_{p}(z/y_{j})}{\Theta_{p}(t^{-1}z/y_{j})}\right]_{1}.
\end{gather*}
For the elliptic Macdonald operator and the elliptic kernel function $\Pi_{MN}(q,t,p)(x,y)$, we have the
following functional equation
\begin{gather}
\{H_{M}(q,t,p)_{x}-t^{M-N}H_{N}(q,t,p)_{y}\}\Pi_{MN}(q,t,p)(x,y)
\nonumber
\\
\qquad
=\frac{(1-t^{M-N})(p;p)_{\infty}^{3}}{\Theta_{p}(t)}C_{MN}(p;x,y)\Pi_{MN}(q,t,p)(x,y).
\label{eq3.6}
\end{gather}
\end{Theorem}

\begin{proof}
The proof is straightforward.
Using Theorems~\ref{th2.3} and~\ref{th2.5}, we calculate the matrix element $\langle 0|\phi^{\ast}_{M}(p;x)[\eta(p;z)]_{1}\phi_{N}(p;y)|0 \rangle$ in two dif\/ferent ways as follows
\begin{gather*}
 \langle 0|\phi^{\ast}_{M}(p;x)[\eta(p;z)]_{1}\phi_{N}(p;y)|0 \rangle
=\frac{t^{-M+1}\Theta_{p}(t^{-1})}{(p;p)_{\infty}^{3}}H_{M}(q,t,p)_{x}\Pi_{MN}(q,t,p)(x,y)
\\
\qquad
\phantom{=}
{}+t^{-M}C_{MN}(p;x,y)\Pi_{MN}(q,t,p)(x,y)
\\
\qquad
=\frac{t^{-N+1}\Theta_{p}(t^{-1})}{(p;p)_{\infty}^{3}}H_{N}(q,t,p)_{y}\Pi_{MN}(q,t,p)(x,y)
+t^{-N}C_{MN}(p;x,y)\Pi_{MN}(q,t,p)(x,y).
\end{gather*}
Therefore we obtain Theorem~\ref{theorem3.2}.
\end{proof}

\begin{Remark}
We can check the following:
\begin{gather*}
C_{MN}(p;x,y)=\left[\prod\limits_{i=1}^{M}
\frac{\Theta_{p}(t^{-1}x_{i}z)}{\Theta_{p}(x_{i}z)}\prod\limits_{j=1}^{N}\frac{\Theta_{p}(z/y_{j})}{\Theta_{p}(t^{-1}z/y_{j})}\right]_{1}
\\
\qquad
\xrightarrow[p \to 0]{}
\left[\prod\limits_{i=1}^{M}\frac{1-t^{-1}x_{i}z}{1-x_{i}z}\prod\limits_{j=1}^{N}\frac{1-z/y_{j}}{1-t^{-1}z/y_{j}}\right]_{1}=1.
\end{gather*}
Hence, by taking the limit $p \to 0$, equation~\eqref{eq3.6} reduces to equation~\eqref{eq3.3}.
\end{Remark}

\subsection*{Acknowledgements}

The author would like to thank Koji Hasegawa and Gen Kuroki for helpful discussions and comments.
The author also would like to thank referees for their valuable comments on improvements of the present paper.

\pdfbookmark[1]{References}{ref}
\LastPageEnding

\end{document}